# A possible use of the Kha's protractor

Amelia Carolina Sparavigna
Dipartimento di Fisica,
Politecnico di Torino, Torino, Italy

I have recently proposed that an object, found in an Egyptian tomb and exposed at the Egyptian Museum of Torino, could be a protractor. The tomb was that of architect Kha, supervisor at Deir El-Medina during the 18[th] Dynasty, and his wife Merit. Considering then the object as a protractor, the architect could have used it to measure the angle of the inclined planes.

As discussed in a previous paper [1], in 1906, an Italian archaeological mission found a tomb of the New Kingdom, that survived intact till its discovery. In this tomb, Kha and his wife Merit were buried. Kha was an architect and an important supervisor at Deir El-Medina, for some architectonical projects of three kings of the 18[th] Dynasty [1]. The coffins and objects of the tomb are exposed at the Museo Egizio of Turin [2].

Kha had, for his afterlife, some objects he used during his job. Among his items we see two cubits and some tools for writing. One of the cubits was of wood and could be folded by hinges. This was probably the cubit Kha handled during the job [3]. As told in Ref.4, wood was the material actually used for the cubit rods by architects. Besides these rulers, several other tools were used in the ancient Egypt masonry, such as plumbs, levels and squares [1,3]. As ours, the Egyptian plumbs consisted of a suspended plumb bob [1].

In the Fig.1, it is shown an object from the Kha's Tomb supposed to be the case of a balance scale (see Fig.1), or the scale itself as reported by the corresponding label. The tool has a complex decoration, that suggested me the case could be used as a protractor, to determine directions and measure angles. The detail of decoration is shown in Fig.2: we see the 16-fold symmetry of a compass rose with 16 leaves. Outside this rose there is a polygonal line with 18 corners and then 36 sides. As I noted in [1], the fraction 1/16, corresponding to one leaf of the decoration, is one of the fractions of the Eye Of Horus, defined during the Old Kingdom to represent the number one, and the number 36, the number of sides, corresponds to the number of Decans, the 36 groups of stars which rise in succession from the horizon due to the earth rotation. Considering the decoration of the case as a protractor, it had two scales, one based on Egyptian fractions, the other based on Decans.

Assuming the tool used to evaluate angles, let me show a possible set-up to measure the angle of an inclined plane. Note that the case has a lid (Fig.1). Let us imagine to remove the lid; the case has a perfectly linear side, that can be put on a smooth surface, as shown in Fig.3. When the surface is horizontal, using a plumb to have the vertical direction, one of the directions of the rose of Fig.2 coincides with the direction of the plumb (the line between leaves 1-16 or 8-9 seems to be perpendicular to the side of the case). If the surface is inclined, the direction of the rose is inclined forming a certain angle with respect to the vertical. This angle has the same value of the angle of the inclined plane (Fig.3 and Fig.4).

As written in [4], geometry originated as a practical science, to measure lengths, surfaces and volumes. In my opinion, Kha could have used his tool, with the contemporary use of a plumb to have a reference direction, to create a set-up able to provide a practical measurement of inclination. The Kha's protractor could be one of the first devices to measure angles.

2. Turin Egyptian Museum: the tomb of Kha
3. Building in Egypt; Pharaonic Stone Masonry, Dieter Arnold, Chapter 6, Tools and their applications, New York and Oxford, 1991
4. http://en.wikipedia.org/wiki/Geometry

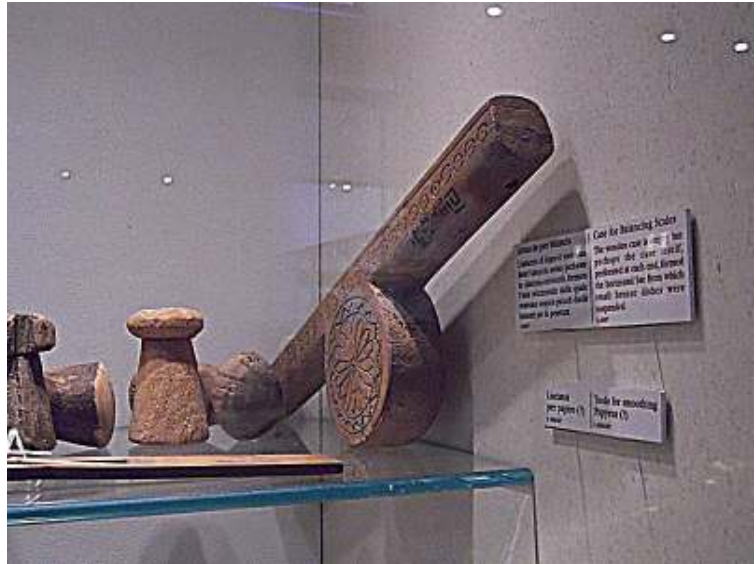
Fig.1 The label is telling that the object is the case of a balance scale (Egyptian Museum, Torino)

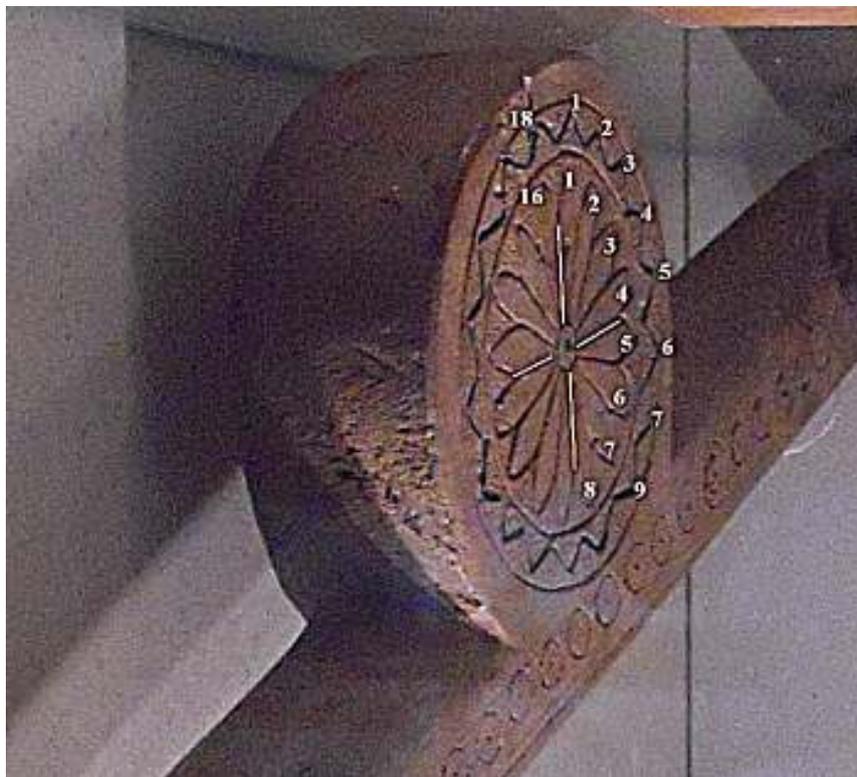
Fig.3 The decoration of this object has a 16-fold symmetry, as a compass rose with 16 leaves. Outside this rose, there is a polygonal line with 18 corners, pointing outwards. They correspond to the same number of corners, pointing inwards. That is, we have a line with 36 sides. The line between leaves 1-16 or 8-9 seems to be perpendicular to the base of the object.

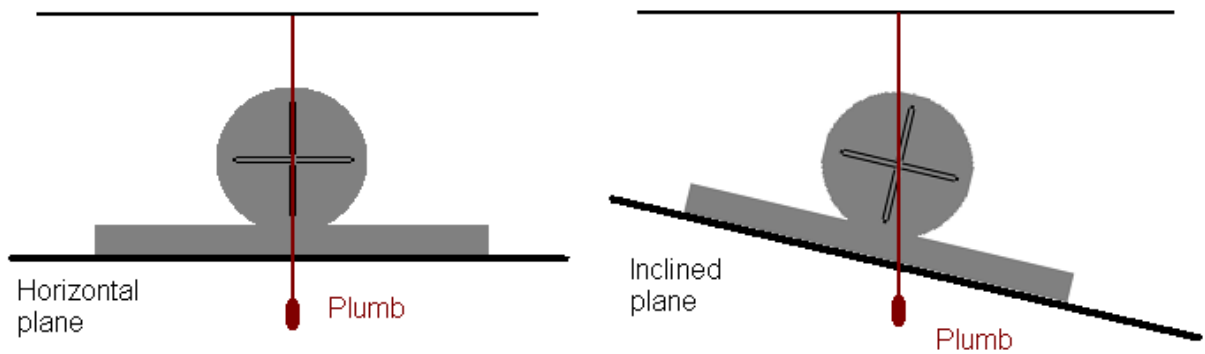

Fig.3. Assuming the object of the Kha's tomb as a protractor, we can measure the angle of an inclined plane. We can put it on a smooth surface. In the case that the surface is horizontal, using a plumb to have the vertical direction, one of the directions of the rose in Fig.2 seems to coincide with the plumb. In the case that the surface is inclined, the direction of the rose is inclined with respect to the vertical.

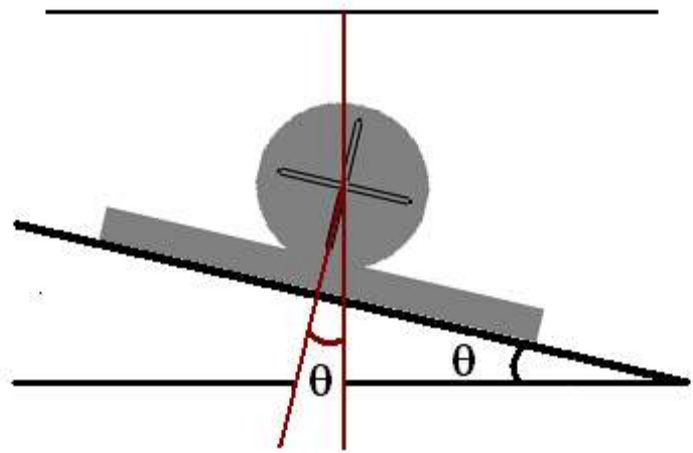

Fig.4: The angle with the vertical direction has the same value of the angle of the inclined plane.